\documentclass[11pt,a4paper]{article}
\usepackage[T1]{fontenc}
\usepackage[utf8]{inputenc}
\usepackage{amsmath,mathrsfs,amssymb,amsfonts,nicefrac}
\usepackage[english,francais]{babel}
\usepackage[top=4cm, bottom=4cm, left=2.5cm, right=2.5cm]{geometry}

\newtheorem{theorem}{Theorem}[section]

\newtheorem{e-proposition}[theorem]{Proposition}

\newtheorem{e-definition}[theorem]{Definition\rm}

\newtheorem{theoreme}{Th\'eor\`eme}[section]
\newtheorem{lemme}[theoreme]{Lemme}

\newcommand{\R}{\mathbb{R}}
\newcommand{\dis}{\displaystyle}
\newcommand{\abs}[1]{\left|#1\right|}
\newcommand{\eps}{\varepsilon}

\renewcommand{\a}{\alpha}
\renewcommand{\leq}{\leqslant}
\renewcommand{\geq}{\geqslant}
\renewcommand{\u}{\underline}
\renewcommand{\o}{\overline}
\renewcommand{\tilde}{\widetilde}

\def\og{\leavevmode\raise.3ex\hbox{$\scriptscriptstyle\langle\!\langle$~}}
\def\fg{\leavevmode\raise.3ex\hbox{~$\!\scriptscriptstyle\,\rangle\!\rangle$}}

\begin{document}
\title{Propagation in Fisher-KPP type equations with fractional diffusion in periodic media}
\author{Xavier CABR\'{E}$^{\small{*}}$, Anne-Charline COULON$^{\small{**}}$, Jean-Michel ROQUEJOFFRE$^{\small{**}}$\\
\footnotesize{$^{\small{*}}$ ICREA and Universitat Politècnica de Catalunya, Dep. de Matemàtica Aplicada I,}\\
\footnotesize{Av. Diagonal 647, 08028 Barcelone, Espagne}\\
\footnotesize{ $^{\small{**}}$Institut de Math\'ematiques, (UMR CNRS 5219), Universit\'e Paul Sabatier,}\\
\footnotesize{118 Route de Narbonne, 31062 Toulouse Cedex, France}}
\date{}
\maketitle

%

\selectlanguage{english}
\begin{abstract}
We are interested in the time asymptotic location of the level sets of solutions to Fisher-KPP reaction-diffusion equations with fractional 
diffusion in periodic media. We show that the speed of propagation is exponential in time, with a precise exponent depending on a periodic 
principal eigenvalue,  and that it does not depend on the space direction. This is in contrast with the Freidlin-G\"artner formula for the standard Laplacian.
\end{abstract}
\vskip 0.5\baselineskip

\selectlanguage{francais}
\begin{abstract}
{\bf Propagation dans les equations de type Fisher-KPP avec diffusion fractionnaire en milieux périodiques. }
On s'int\'eresse ici \`a la localisation asymptotique en temps des lignes de niveaux des solutions d'equations de réaction-diffusion de type 
Fisher-KPP avec diffusion fractionnaire en milieu périodique. Nous montrons que la vitesse de propagation est exponentielle en temps, avec un
 exposant dépendant d'une valeur propre principale périodique, et que cette vitesse ne d\'epend pas de la direction de propagation. 
Ceci est en contraste avec la formule de Freidlin-G\"artner pour le laplacien standard.
\end{abstract}

\selectlanguage{francais}
\section*{Version fran\c{c}aise abr\'eg\'ee}
Consid\'erons l'equation de réaction-diffusion de type Fisher-KPP (Kolmogorov - Petrovskii - Piskunov) suivante:
\begin{equation}\label{equa1}
           u_t+(-\Delta)^{\a}u= \mu(x) u-u^2, \ \ \ x\in \R^d, t>0,
\end{equation}
avec donn\'ee de Cauchy  $u(\cdot,0)=u_0$. Dans ce modèle, $(-\Delta)^{\a}$ repr\'esente le laplacien fractionnaire d'ordre $\a \in (0,1);$ 
 la fonction $\mu$ est supposée périodique en chaque variable $x_i$ d'espace et vérifie $0<\min \mu \leq \mu(x)$. Les conditions sur la donnée 
initiale $u_0$ sont données dans le théorème principal. Soit $\lambda_1$ la valeur propre principale p\'eriodique de $(-\Delta)^\alpha-\mu(x)I$. 
Dans le cas $\lambda_1 \geq 0$, la solution tend vers 0 quand $t \rightarrow +\infty$ (voir~\cite{BRR}); c'est pourquoi on va supposer $\lambda_1<0$. 
Il existe alors une unique solution stationnaire strictement positive pour \eqref{equa1}, not\' ee  $u_+$. Par unicité, $u_+$ est périodique. Si $u$ est solution de 
\eqref{equa1} alors $u(x,t)\to u_+(x)$ quand $t\to+\infty$, uniform\'ement sur tout compact et on cherche \`a comprendre \`a quelle vitesse l'\'etat  stable $u_+$ 
envahit l'\'etat instable 0. 

Lorsque $\a=1$, beaucoup de résultats sont connus. En milieu homogène (i.e. $\mu$ est constante, par exemple $\mu \equiv 1$ et donc $u_+\equiv1$), 
Aronson et Weinberger \cite{AW} prouvent que les lignes de niveau d'une solution $u$ issue d'une donn\'ee \`a support compact se propagent asymptotiquement 
en temps à vistesse constante \'egale \`a $2$, indépendamment de la direction de propagation.
Quand $u$ est issue d'une donn\'ee faiblement d\'ecroissante $u_0$, la vitesse asymptotique est exponentielle en temps (voir \cite{HR}).
 En milieu périodique (i.e. $\mu$ est non constante et périodique), si $u_0$ est à support compact, Freidlin et
Gärtner \cite{FG} prouvent, avec des outils probabilistes, qu'il existe une vitesse notée $w^*(e)$ dans chaque direction de propagation 
$e \in S^{d-1}$ (la d\'efinition et  l'expression de $w^*(e)$ sont données en \eqref{GFformula}). 
D'autres d\'emonstrations, utilisant des arguments de syst\`emes dynamiques ou d'EDP, sont propos\'ees dans \cite{Wein}, \cite{ES1}, \cite{BHN}.

Lorsque $\a \in (0,1)$  et dans le cas homogène $(\mu\equiv 1)$, il est prouvé dans \cite{JMRXC}, \cite{JMRXCbis} que les lignes de niveaux de $u$ 
se propagent approximativement \`a la vitesse $e^{\frac{t}{d+2\a}}$. La transition entre les propagations linéaire ($\a =1$) et exponentielle ($\a \in (0,1)$) 
est traitée dans \cite{CR}. Ici nous généralisons ces résultats en milieu périodique. Contrairement au cas du laplacien standard, la vitesse de propagation ne dépend
pas de la direction:\\
\begin{theoreme}
Supposons $\lambda_1 <0$ et considérons $u$ la solution de \eqref{equa1} où la donnée initiale $u_0$ est continue par morceaux, positive,
$u_0 \not \equiv 0$  et $u_0(x)={\rm{O}}( \abs{x}^{-(d+2\a)})$ quand $\abs{x}\to +\infty$.
Alors pour tout $\lambda \in (0,\min \mu)$, il existe une constante $c_{\lambda}>0$ et un temps $t_{\lambda}>0$ (ces deux constantes dépendent de $\lambda$ et $u_0$)
 tels que pour tout $t \geq t_{\lambda}$:
$$ \{ x\in \R^d \ | \ u(x,t) = \lambda\} \subset \{x \in \R^d \ | \  c_{\lambda}e^{\frac{{ \abs{\lambda_1}} }{d+2\a}t} \leq \abs{x} \leq 
c_{\lambda}^{-1}e^{\frac{{\abs{\lambda_1}} }{d+2\a}t}\}.$$
\end{theoreme}
La preuve est basée sur la construction de sous-solutions et sur-solutions explicites. Dans la version anglaise de cette note, 
nous pr\'esentons une d\'emonstration compl\`ete du cas $\a<\frac{1}{2}$ (le cas g\'en\'eral sera trait\'e dans \cite{CCR}). 
L'étape principale est le lemme suivant: \\
\begin{lemme}

Supposons $\lambda_1<0$ et $\a < \frac{1}{2}$. Pour des constantes strictement positives $\u a$, $\o a$, $\u B$, $\o B$ et $M$, soit
 $$\u{u}(x,t)=\dis{\frac{\u a \phi_1(x)}{\abs{\lambda_1}^{-1}+ \u b(t)\abs{x}^{d+2\a}}} \mbox{ \ \ et \ \ }\o{u}(x,t)=\dis{\frac{\o a \phi_1(x)}
{\abs{\lambda_1}^{-1}+ \o b(t)\abs{x}^{d+2\a}}},$$ avec $\u b(t)=(M\abs{\lambda_1}^{-1}+\u B^{-\frac{2\a}{d+2\a}}e^{\frac{2\a\abs{\lambda_1}}{d+2\a}t})^{-\frac{d+2\a}{2\a}}$
 et $\o b(t)=(-M\abs{\lambda_1}^{-1}+\o B^{-\frac{2\a}{d+2\a}}e^{\frac{2\a\abs{\lambda_1}}{d+2\a}t})^{-\frac{d+2\a}{2\a}}$. \\
Il existe une constante $M>0$ telle que:\\
$\bullet$ si $0<\u B< (\abs{\lambda_1}M^{-1})^{\frac{d+2\a}{2\a}}$ et $0<\u a\leq (\max \phi_1)^{-1}( 1-\u B^{\frac{2\a}{d+2\a}}M\abs{\lambda_1}^{-1})$,
 alors $\u{u}$ est sous-solution de \eqref{equa1} pour $t>0$.\\
$\bullet$ pour toute constante $\o B >0$, si $t_0=(d+2\a)(2\a \abs{\lambda_1})^{-
1}\ln(2^{-1}+M\abs{\lambda_1}^{-1}\o B^{\frac{2\a}{d+2\a}})$ et $\o a \geq (\min \phi_1)^{-1}(1+2\o B^{\frac{2\a}{d+2\a}}M\abs{\lambda_1}^{-1})$, 
alors $\o{u}$ est une sur-solution de \eqref{equa1} pour $t > t_0$.

\end{lemme}

Un lemme analogue est vrai dans le cas $\a \geq \frac{1}{2}$. L'id\'ee sous-jacente \`a la construction de  $\u{u}$ et $\o{u}$ est expos\'ee dans la 
version anglaise de cette note.

\selectlanguage{english}
\section{Introduction: Motivation and main result}
We are interested in the time asymptotic location of the level sets of solutions to the equation
\begin{equation}\label{equa}
           u_t+(-\Delta)^{\a}u= \mu(x) u-u^2,\ \ \ x\in \R^d, t>0,
\end{equation}
with initial condition  $u(\cdot,0)=u_0,$ where $\a \in (0,1)$, $\mu$ is periodic in each $x_i$-variable and satisfies  $0<\min \mu \leq \mu(x)$, 
and  $(-\Delta)^\alpha$ is the fractional Laplacian.
The nonlinearity $\mu(x)u-u^2$ is often referred to as a Fisher-KPP type nonlinearity.

Let ${\lambda_1}$ be the principal periodic eigenvalue of the operator $(-\Delta)^{\a} -\mu(x)I$.  From \cite{BRR} we know that if ${\lambda_1}\geq 0$, 
every solution to \eqref{equa} starting with a bounded nonnegative initial condition tends to $0$ as $t\to+\infty$. Thus in this paper we assume ${\lambda_1}<0$. 
Then, by \cite{BRR}, the solution to \eqref{equa} tends, as $t\to+\infty$, to  the unique bounded positive steady solution to \eqref{equa}, denoted by $u_+$. 
By uniqueness, $u_+$ is periodic. The convergence holds on every compact set. Hence, the level sets of $u$ spread to infinity for large times, and we wish to 
understand how fast. To do it, we look for a function $R_e(t)$ going to $+\infty$ as $t$ tends to $+\infty$ such that, for every direction $e\in{S}^{d-1}$ and 
every constant $c\in (0,1)$,
 \begin{equation}\label{Re}
 \liminf_{t\to+\infty}\left(\inf_{\{ x= \rho e, 0\leq  \rho \leq  R_e(ct)\}}u(x,t)\right)>0 \  \hbox{ and }\  \  \limsup_{t\to+\infty}\left(\sup_{\{ x= \rho e, 
\rho \geq R_e(c^{-1}t)\}}u(x,t)\right) =0.
\end{equation}
We adopt this slightly unusual definition to cover both linear and exponential propagation.

The case $\a=1$ is well studied. In homogeneous media, when the function $\mu$ is constant (say equal to~$1$), \cite{AW} establishes that, if $u_0$
is compactly supported, then we may choose $R_e(t)=2t$ regardless of the direction $e$ of propagation. See  \cite{HR} for the case of slowly decreasing 
initial conditions $u_0$ (here, essentially, any sufficiently rapidly increasing function becomes an $R_e(t)$ for some $u_0$).
In space periodic media, starting from a compactly supported initial data, Freidlin and Gärtner
\cite{FG} have characterised $R_e(t)$ by
\begin{equation}
\label{GFformula}
R_e(t)=w^*(e)t, \  \    \  w^*(e)=\min_{e' \in S^{\tiny{d-1}},e'\cdot e>0}\frac{c^*(e')}{e'\cdot e},
\end{equation}
where $c^*(e')$ is the minimal speed of pulsating travelling fronts in the direction $e'$. Their proof uses probabilistic tools; 
proofs using dynamical systems or PDE arguments are given in  \cite{Wein}, \cite{ES1}, \cite{BHN}.

For $\alpha\in(0,1)$ and $\mu \equiv 1$ in \eqref{equa}, propagation is exponential in time. Although this fact was well noted in physics references, 
the first mathematically rigorous result is \cite{JMRXC}, \cite{JMRXCbis}, which proved the exponential propagation in time; in fact, that $R_e(t)=e^{\frac{t}{d+2\a}}$  
in this case. The transition between linear propagation ($\a=1$) and exponential propagation ($\a\in (0,1)$) is examined in \cite{CR}.  See \cite{Jimmy} for equations 
of the type \eqref{equa} with nonsingular integral dispersal. 

 In this note we prove the following result :
 \begin{theorem}\label{thm1}
 Assume that $\lambda_1 <0$. Let $u$ be the solution to \eqref{equa} with $u_0$ piecewise continuous, nonnegative,
$u_0 \not\equiv 0$, and $u_0(x)= {\rm{O}} (\abs{x} ^{-(d+2\a)})$ as $\abs x \rightarrow \infty$.
Then, for every $\lambda \in (0, \min \mu)$, there exist $c_{\lambda}>0$  and a time $t_{\lambda}>0$ (all depending on $\lambda$ and $u_0$) such that, 
for all $t \geq t_{\lambda}$,
\begin{equation}\label{ineg}
 \{ x\in \R^d \ | \ u(x,t) = \lambda\} \subset \{x \in \R^d \ | \  c_{\lambda}e^{\frac{{ \abs{\lambda_1}} }{d+2\a}t} \leq \abs{x} \leq c_{\lambda}^{-1}
e^{\frac{{\abs{\lambda_1}} }{d+2\a}t}\}.
\end{equation}
\end{theorem}
Theorem \ref{thm1} gives that \eqref{Re} holds with $R_e(t)=e^{\frac{\abs{\lambda_1} }{d+2\a}t}$. Thus, spreading does not depend on the direction of propagation,
 and this is in contrast with \eqref{GFformula} for the standard Laplacian.  Moreover, the estimate  that we obtain is much sharper than that in \cite{JMRXC}, 
\cite{JMRXCbis} for $\mu \equiv 1$. Indeed, to guarantee the limits in \eqref{Re}, \cite{JMRXCbis} needed to assume $\abs x \leq Ce^{\sigma_1t}$ 
(respectively $\abs x \geq Ce^{\sigma_2t}$) with $\sigma_1<\frac{1}{d+2\a}<\sigma_2$. Note that $\lambda_1=-1$ when $\mu\equiv 1$. When 
$\lambda \geq \min u_+$, \eqref{ineg} can not hold since $u(x,t) \rightarrow u_+(x)$ as $t \rightarrow +\infty$. On the other hand one may expect it to hold for 
$\lambda \in (0, \min u_+)$. Here we prove it for $\lambda \in (0,\min \mu)$; it is easy to see that  $\min \mu \leq \min u_+$.

The proof of Theorem \ref{thm1} is quite simple: it relies on the construction of explicit subsolutions and supersolutions, which are themselves based on a 
nonlinear transport equation, \eqref{e2.1}, satisfied asymptotically by a correctly rescaled version of the solution $u$. The rest of this note is devoted 
to a full proof in the case $\a< \frac{1}{2}$ but the result remain true for all $\a \in (0,1)$ and this will be explained in \cite{CCR}.

\section{The proof of Theorem \ref{thm1}}
We will from now on assume $\a< \frac{1}{2}$. Recall that ${\lambda_1}<0$ denotes the principal periodic eigenvalue of the operator $(-\Delta)^{\a} -\mu(x)I$ 
and that the corresponding periodic eigenfunction is denoted by $\phi_1$. \\
\textit{Step 1.} Let us write $u(x,t)=\phi_1(x) v(x,t)$ and define $w(y,t)=v(yr(t),t)$ for $r(t)=e^{\frac{\abs{\lambda_1}t}{d+2\a}}$.
Thus, for $y\in \R^d$ and $t>0$, $w$ solves
$$
 w_t- \frac{\abs{\lambda_1}}{d+2\a} y\cdot w_y+e^{\frac{-2\a\abs{\lambda_1}t}{d+2\a}}\left[(-\Delta)^{\a}w-\frac{Kw}{ \phi_1(yr(t))}\right]=\abs{\lambda_1} 
w - \phi_1(yr(t)) w^2, 
$$
where we have used $\lambda_1<0$ and we have defined
$$Kw(y)= C_{d,\a} PV \dis{\int_{\R^d}} \dis{\frac{ \phi_1(yr(t))- \phi_1 (\o yr(t))}{\abs{y-\o y}^{d+2\a}}} (w(y)-w(\o y))d\o y.
$$
If we formally neglect the term 
$e^{\frac{-2\a\abs{\lambda_1}t}{d+2\a}}\left[(-\Delta)^{\a}w-\frac{Kw}{ \phi_1(yr(t))}\right]$ which should go to $0$ as $t \rightarrow +\infty$,
 we get the transport equation
\begin{equation}
\label{e2.1}
\tilde{w}_t- \frac{\abs{\lambda_1}}{d+2\a} y\cdot {\tilde{w}}_y=\abs{\lambda_1}{\tilde{w}} - \phi_1(yr(t)) {\tilde{w}}^2, \ \ y\in \R^d, t>0.
\end{equation}
Equation \eqref{e2.1}, completed by an initial datum $\tilde{w_0}$, is solved as:
$$\tilde{w}(y,t)=\frac{1}{  \phi_1(yr(t))\abs{\lambda_1}^{-1}+e^{-\abs{\lambda_1}t}\left(\tilde{w_0}(yr(t))^{-1}-  \phi_1(yr(t))\abs{\lambda_1}^{-1}\right)}.
$$
Taking into account (see for instance \cite{JMRXCbis}) that $\abs{x}^{d+2\a} u(x,t)$ is uniformly bounded from above and below (but of course not uniformly in $t$),
  it is natural to specialise
$\tilde{w_0}(y)=\frac{1}{1+\abs{y}^{d+2\a}}$. In this case we have 
$$\tilde{w}(y,t)=\frac{1}{ \phi_1(yr(t))\abs{\lambda_1}^{-1} (1-e^{-\abs{\lambda_1}t})+e^{-\abs{\lambda_1}t}+\abs{y}^{d+2\a}}.
$$
Since $\phi_1$ is bounded above and below and $t$ tends to $+\infty$, coming back to the function $v(x,t)=w(x r(t)^{-1},t)$, the idea is to consider the following 
family of functions modelled by $\tilde{w}$:
\begin{equation}
\label{e2.2}
\tilde{v}(x,t)=\frac{a}{\abs{\lambda_1}^{-1}+b(t)\abs{x}^{d+2\a}}, \qquad \ \tilde{u}(x,t)=\phi_1(x)\tilde{v}(x,t),
\end{equation}
and to adjust  $a>0$ and $b(t)$ asymptotically proportional to $e^{-\abs{\lambda_1}t}$ so that the function $\tilde{u}(x,t)$ serves as a subsolution or a supersolution 
to \eqref{equa}.\\
\textit{Step 2.} Thus, let $\tilde{v}(x,t)$ be defined by \eqref{e2.2}. 
Let us consider the operator $ \tilde{K}$ defined by \\ $$ \tilde{K}\tilde g(x)= C_{d,\a} PV\dis{ \int_{\R^d} \frac{\phi_1(x)-\phi_1(\o x)}{\abs{x-\o x}^{d+2\a}} 
(\tilde g(x)-\tilde g(\o x))d\o x}.$$
To prove the result we need to understand the effect of the operators $(-\Delta)^{\a}$ and $\tilde{K}$ on $\tilde{v}$. From simple computations there is a constant 
$D>0$ such that 
\begin{equation}\label{estimation}
 \vert(-\Delta)^{\a} \tilde{v}\vert \leq Db(t)^{\frac{2\a}{d+2\a}}\tilde{v} \ \mbox{ and }  \ \vert \tilde{K}\tilde{v}\vert \leq Db(t)^{\frac{2\a}{d+2\a}}\tilde{v} 
\quad \mbox{ in } \R^d. 
 \end{equation}
Note that the estimate on $\tilde{K}$ does not hold for $\a \geq \frac{1}{2}$ and we must replace the term 
$b(t)^{\frac{2\a}{d+2\a}}$ by $b(t)^{\frac{2\a-\gamma}{d+2\a}}$, where $\gamma$ is any positive constant between $2\a-1$ and $1$.

We use these estimates to find $a$ and $b(t)$, that we denote by $\u a$ and $\u b(t)$, such that $\u u(x,t)= \phi_1(x)\tilde{v}(x,t) $ is a subsolution to \eqref{equa}. 
Using $\lambda_1<0$, note that
\begin{eqnarray}{ll}
\u u_t +(-\Delta)^{\a}\u u -\mu(x) \u u +\u u^2 &=& \phi_1 \tilde{v}_t +  \phi_1(-\Delta)^{\a}\tilde{v} - 
\tilde{ K}\tilde{v}-\abs{\lambda_1}  \phi_1 \tilde{v}+ \phi_1^2 \tilde{v}^2   \nonumber\\
&  \leq&  \phi_1\tilde{v}^2 \u a^{-1}\left\{ -{\u b'}(t) \abs{x}^{d+2\a} 
+ D\left((\min \phi_1)^{-1}+1\right)\u b(t)^{\frac{2\a}{d+2\a}}\right.\\ \label{sousol}
&&\left.\left(\abs{\lambda_1}^{-1}  +\u b(t)\abs{x}^{d+2\a}\right)-\abs{\lambda_1} 
\left(\abs{\lambda_1}^{-1}+\u b(t)\abs{x}^{d+2\a}\right)   + \u a \phi_1\right\}. \nonumber
\end{eqnarray}
We define $M=D\left((\min \phi_1)^{-1}+1\right)>0$ and we take a constant $\u B> 0$. The solution to
\begin{equation}
\label{e2.3}
-{\u b'}(t)+ M\u b(t)^{\frac{2\a}{d+2\a}+1}-\abs{\lambda_1}\u b(t)=0,\   \    \    \    \u b(0)=\left(M\abs{\lambda_1}^{-1}+
\u B^{-\frac{2\a}{d+2\a}}\right)^{-\frac{d+2\a}{2\a}},
\end{equation}
is $\u b(t)=\left(M\abs{\lambda_1}^{-1}+\u B^{-\frac{2\a}{d+2\a}}e^{\frac{2\a\abs{\lambda_1}}{d+2\a}t}\right)^{-\frac{d+2\a}{2\a}}$ and we have 
$\u b(t) \leq \u Be^{-\abs{\lambda_1} t} \leq \u B$ for all $t \geq 0$. Thus if we choose $\u B>0$ and $\u a>0$ such that 
$\u B< (M^{-1}\abs{\lambda_1})^{\frac{d+2\a}{2\a}}$ and $ \u a\leq (\max \phi_1)^{-1}( 1-\u B^{\frac{2\a}{d+2\a}}M\abs{\lambda_1}^{-1}),
$
the right hand side of inequality \eqref{sousol} is less than or equal to $0$.

We make $\o u(x,t)=\phi_1(x)\tilde{v}(x,t) $, with $a$ and $b(t)$ denoted by $\o a$ and $\o b(t)$, a supersolution in a similar fashion. 
Here we must replace $D$ and $M$ by $-D$ and $-M$ in \eqref{sousol} and \eqref{e2.3}. Given any positive constant $\o B>0$, we take
 $\o b(t)=(-M\abs{\lambda_1}^{-1}+\o B^{-\frac{2\a}{d+2\a}}e^{\frac{2\a\abs{\lambda_1}}{d+2\a}t})^{-\frac{d+2\a}{2\a}}$  with initial datum  
$\o b(0)=(-M\abs{\lambda_1}^{-1}+\o B^{-\frac{2\a}{d+2\a}})^{-\frac{d+2\a}{2\a}}$. Defining $t_0= (d+2\a)(2\a \abs{\lambda_1})^{-1}\ln(2^{-1}+
M \abs{\lambda_1}^{-1}\o B^{\frac{2\a}{d+2\a}}),$ we have $\o B e^{-\abs{\lambda_1}t} \leq \o b(t) \leq 2^{\frac{d+2\a}{2\a}}\o B$ for all $t \geq t_0$. 
Thus, taking $\o a \geq (\min \phi_1)^{-1}(1+2\o{B}^{\frac{2\a}{d+2\a}}M\abs{\lambda_1}^{-1}),$
we get a supersolution to \eqref{equa} for $t>t_0$.\\
\textit{Step 3.} Let us prove the main theorem in the case $\a< \frac{1}{2}$. Due to the assumption on $u_0$ for large values of $x$, we can choose $\o a$
and $\o B$ satisfying the conditions obtained in step 2 and such that $\o{u}(x,t_0) \geq u_0(x)$. Thus, given any $\lambda>0$ and using the maximum principle, 
we get for $t>0$
$$\{x \in \R^d \ | \ \abs{x} >C_2e^{\frac{\abs{\lambda_1} }{d+2\a}t}\} \subset \{ x\in \R^d \ | \ u(x,t) <\lambda\},$$
where $C_2^{d+2\a}=\lambda^{-1} \o B^{-1} e^{\abs{\lambda_1}t_0}\o a \max \phi_1$.

The last item to prove is that, for all $\lambda \in (0, \min \mu)$, there exist $c_1>0$ and a time $t_1>0$ (both depending on $\lambda$) such that for $t\geq t_1$
\begin{equation}\label{incl}
\{x \in \R^d \ | \ \abs{x} < c_1e^{\frac{\abs{\lambda_1} }{d+2\a}t}\} \subset \{ x\in \R^d \ | \ u(x,t) >\lambda\}.
 \end{equation} 
Since $\u{u}(\cdot,0) \leq u_0$ may not hold, we look for a time $t_1>\max((\min u_+)^{-\frac{2\a}{d}},2M\abs{\lambda_1}^{1+\frac{2\a}{d+2\a}})$ such that 
$u(x,t_1) \geq 2^{-1} \min u_+$ for $\abs x \leq 1$. Moreover, we know that
$u(x, t_1) \geq \frac{ct_1}{t_1^{\frac{d}{2\a}+1}+ \abs{x}^{d+2\a}}$ for $\abs{x} \geq 1$, where $c\in (0,1)$ is a constant (see the proof of Lemma 2.2 in 
\cite{JMRXCbis} for the computation). Consequently we choose $\u a=\frac{c}{2 t_1^{\frac{d}{2\a}}\abs{\lambda_1}\max \phi_1 }$ and 
$\u B = \frac{2^{\frac{d+2\a}{2\a}}}{\abs{\lambda_1}t_1^{\frac{d}{2\a}+1}}$ to have $\u{u}(x,0) \leq u(x,t_1)$ for all $x \in \R^d$. 
Taking $t_1$ larger if necessary, the requirements for $\u B$ and $\u a$ in step 2 are satisfied.
By the maximum principle $\u {u}(x,t-t_1)\leq u(x,t)$ for $t \geq t_1.$
Let us define $\eps=\frac{c \min \phi_1 }{4t_1^{\frac{d}{2\a}} \max \phi_1}$ and $\tilde{c_1}^{d+2\a}= e^{-\abs{\lambda_1} t_1}\abs{\lambda_1}^{-1} \u B^{-1}$, 
and take $x$ such that $\abs{x} < \tilde{c_1}e^{\frac{\abs{\lambda_1} }{d+2\a}t}$. Since $\u b(t) \leq \u B e^{-\abs{\lambda_1}t}$,  for $t \geq t_1$ we get
$$u(x,t)> \frac{\u a \min \phi_1}{\abs{\lambda_1}^{-1}+\u Be^{\abs{\lambda_1}t_1}\tilde{c_1}^{d+2\a}}= \frac{\u a\min \phi_1}{2\abs{\lambda_1}^{-1}} = \eps.
$$

Thus, we have found an $\eps >0$ such that \eqref{incl} holds for $t>t_1$ (with $c_1$ and $\lambda$ replaced by $\tilde{c_1}$ and $\eps$ respectively). 
Note that $\mu(x) u-u^2\geq (\min \mu)u-u^2 = u(\min \mu -u)$. Now, we use Lemma 3.3 in \cite{JMRXCbis} (which also holds for supersolutions), replacing 
the stable state $1$ by $\min \mu$. For every $\lambda\in(0,\min \mu)$, we deduce the existence of a positive constant $c_1$ and a time $t_{\lambda}$ 
(both depending on $\lambda$) such that $u(x,t)> \lambda$ for $ \abs{x} < c_1e^{\frac{\abs{\lambda_1} }{d+2\a}t} $ and $t\geq t_{\lambda}$ .
Theorem \ref{thm1} is proved by taking $c_{\lambda}=\min(c_1,C_2^{-1})$.
\section{Concluding remarks}
Clearly, the transport equation \eqref{e2.1} is a time attractor to the renormalised solution $w$. In the forthcoming paper \cite{CCR} we will study this question 
in a more precise way.

The nonlinearity $\mu(x)u-u^2$ is quite special, it yields simple solutions to \eqref{e2.1}. For more general nonlinearities $f(x,u)$, it is more intricate to find 
sub and supersolutions of the form given by $\tilde{u}$ in \eqref{e2.2}, and the proof is less transparent than in the present case. This will be developped in 
\cite{CCR}, as well as more general diffusions than the fractional diffusion.

\section*{Acknowledgements}
X. Cabr\'e was supported by the Spain Research projects MTM2008-06349-C03-01, MTM2011-27739-C04-01
and the Catalonia Reasearch project 2009SGR345. A.-C. Coulon and J.-M. Roquejoffre were supported by ANR grant PREFERED.

\end{document}